\documentclass{amsart}
\pdfoutput=1

\usepackage{amsmath}

\usepackage{url}

\usepackage{hyperref}

\begin{document}

\title[Letters from Burnside to Fricke and the Burnside Problem]{Letters from William Burnside to Robert Fricke:
Automorphic Functions, and the Emergence of the Burnside Problem}
\author{Clemens Adelmann}
\author{Eberhard H.-A. Gerbracht}
\address{Technische Universit\"at Braunschweig, Institut f\"ur Analysis und Algebra, Pockels\-stra\ss e 14, 38106 Braunschweig, Germany}
\email{c.adelmann@tu-bs.de}
\email{e.gerbracht@web.de}

\date{\today}


\begin{abstract}
Two letters from William Burnside have recently been found in the
{\em Nachlass} of Robert Fricke that contain instances of
Burnside's Problem prior to its first publication. We present
these letters as a whole to the public for the first time. We
draw a picture of these two mathematicians and describe their
 activities leading to their correspondence. We thus gain an
insight into their respective motivations, reactions, and
attitudes, which may sharpen the current understanding of
professional and social interactions of the mathematical
community at the turn of the 20th century.
\end{abstract}

\maketitle


\section{The simple group of order 504 -- a first meeting of  minds}
Until 1902, when the publication list of the then fifty-year-old
William Burnside already encompassed 90 papers, only one of these
had appeared in a non-British journal: a three-page article
\cite{Burnside504} entitled ``Note on the simple group of order
504'' in volume 52 of {\em Mathematische Annalen} in 1898. In
this paper Burnside deduced a presentation of that group in terms
of generators and relations,\footnote{Burnside's presentation is
given by two generators $A, B$ and relations
$$A^7=B^2=(AB)^3=(A^3BA^5BA^3B)^2=E,$$
where $E$ denotes the unit element in the group.} which is based
on the fact that it is isomorphic to the group of linear
fractional transformations with coefficients in the finite field
with 8 elements. The proof presented is very concise and terse,
consisting only of a succession of algebraic identities, while
calculations in the concrete group of transformations are omitted.

In the very same volume of {\em Mathematische Annalen}, only one
issue later, there can be found a paper \cite{Fricke504} by Robert
Fricke entitled ``Ueber eine einfache Gruppe von 504
Operationen'' (On a simple group of 504 operations), which is on
exactly the same subject. But though the subject is the same and
the  line of thought presented ultimately leads to the identical
result, the flavour and style of exposition are completely
different: Fricke's paper consists of 19 pages, and is
accompanied by six very detailed pictures, most of which are
larger than half a printed page. Here the author concentrated on
giving a geometric description of the group in question as a
discrete group acting on the hyperbolic plane. Most of the
necessary calculations are read off the illustrations.

Two mathematicians publishing in the same journal on the same
subject -- per se, this would not justify the need for any
special attention. The different approaches observed may be seen
as a testament to the fact that one of the authors has used
``superior'' methods to simplify and shorten his line of
arguments. This opinion may be further supported by taking into
consideration that William Burnside today is known by many, and
not by only experts, as ``the first to develop the theory of
groups from a modern abstract point of view'', and a ``giant of
the subject'',\footnote{See \cite{MacTutorBurnside} and
\cite{LamFiniteII}, p.~465, respectively.} and thus has become
both a proponent of mathematical modernity and a figure of nearly
mythic proportions. Robert Fricke, on the other hand, who was as
industrious as Burnside and was held in high esteem by his
contemporaries, was for the most part considered no more than a
faithful disciple of his academic teacher, Felix Klein. Thus
Fricke's achievements were not only mostly overshadowed by the
reputation of Klein, but large parts of his work were regarded as
old-fashioned even by the end of his life.\footnote{See e.g.\ the
opinion of K.~O.~Friedrichs (1901--1982), Fricke's successor at the
{\em Technische Hochschule} in Braunschweig, as described in
\cite{RoweOldGuard}, p.~75.} Consequently these two articles
might be seen as nothing more than a display of the difference
between two mathematicians, style and content giving evidence for
Burnside's ``modernity'' and Fricke's ``old-fashionedness''.

Still, there is more to this chance meeting of two at first sight
so dissimilar mathematicians: the way in which their two articles
complement each other so well, without the authors having met
before, indicates that Burnside and Fricke might have had more in
common than is visible at a cursory glance. If only these
articles are taken as evidence, this statement may seem to be a
little far-fetched. Nevertheless, in the sequel we will report on
another, up until now unknown, but significant occasion, when the
paths of Burnside and Fricke crossed again: in 1900 and 1901 the
two mathematicians exchanged several letters. Two of Burnside's
letters have recently been discovered by the authors in Fricke's
{\em Nachlass}.

The mere existence of this correspondence is in itself remarkable,
since Burnside's workstyle was characterized by a later
biographer as having been conducted with no apparent ``extensive
direct contacts with other mathematicians interested in the
subject [of group theory]'',  to the point that Burnside worked
``in isolation, possibly even more so than was normal for his
times, with little opportunity (or, perhaps, inclination) to
discuss his ideas with others''.\footnote{See
\cite{NeumannBurnside}, p.~32.} Burnside's letters prompt us to
now see his isolationist stance (if it existed at all) in a
slightly different light: it was not universal, but seems to have
been only with respect to certain mathematicians, as e.g.\
Frobenius.\footnote{See \cite{LamFiniteII}, pp.~472ff, and
\cite{NeumannBurnside}, p.~31.} Others, like Fricke, were able to
draw strong reactions from Burnside.

Furthermore, and maybe even more noteworthy, these letters contain the first known instance
of Burnside formulating what was later to be called the {\em Burnside Problem},
and show the state of his research in this regard several months before the publication
of the actual article \cite{BurnsideProblem}.

In our paper, besides presenting these letters as a whole to the
public for the first time, we will try to explain why Burnside
might have put his trust in Robert Fricke in such an
uncharacteristic way. At the same time we will draw a picture of
these  two  quite similar 19th century mathematicians on the
brink of mathematical modernity, and demonstrate their common
mathematical interests in the, at that time still interconnected,
fields of discrete groups and automorphic functions, as well as
their attitudes and their interactions with the mathematical
community.

\section{Robert Fricke, Felix Klein, and automorphic functions}

Since it was Fricke who seems to have initiated the contact with Burnside,
we start with him by giving a short description of his life and
work.\footnote{A more detailed biographical account can be found in \cite{Gerke}.}

Born in Helmstedt, Germany, on 24 September 1861, as the second
of four children of a civil servant, Robert Fricke grew up in
Braunschweig, where he finished school in summer 1880. He then
studied mathematics, physics and philosophy for a teacher's degree
at the universities in G\"ottingen, Z\"urich, Berlin, Stra\ss burg,
and finally Leipzig, where he moved in winter 1883. There he
attended the lectures and seminars of Felix Klein (1849--1925),
and was immediately attracted by Klein's approach to mathematics.
Klein on the other hand soon recognized Fricke's talent, and took
him under his wing. In 1885 Fricke received both a teacher's
degree and a {\em Dr.\,phil.}\ from the University of Leipzig.
He returned to Braunschweig, where he first became a school
teacher, but then was granted leave to work as a private teacher
for the sons of the Prince Regent of the Duchy of Braunschweig.
As a result, having more time to advance his mathematical studies,
he kept in touch with Klein, who by 1886 had accepted a chair in
G\"ottingen. In 1891 Fricke decided to give up school teaching
completely in favour of an academic career. He habilitated at the
university of Kiel in the same year, and there became a {\em
Privatdozent}. In 1892 he moved to G\"ottingen to work nearer to
Klein. In 1894 he was appointed professor of higher mathematics
at the {\em Technische Hochschule} in Braunschweig, where he
succeeded Richard Dedekind (1831--1916). He held this position
until his death in 1930. The bond between Fricke and Klein, which
was furthered by Fricke marrying a niece of Klein's in 1894,
remained strong during their lifetimes.\footnote{More than 400
letters between them are kept in the archives in Braunschweig and
G\"ottingen.}

\medskip
To understand the development leading to the correspondence
between Fricke and Burnside we have to go back in time to the
early 1880s, when Klein had just been appointed professor of
geometry in Leipzig. At that time, in 1881 and 1882, Klein was in
scientific competition with the French mathematician Henri
Poincar\'e (1854--1912). The common subject of their research was
the theory of what would later\footnote{Klein coined the term
{\em automorphic function} about 1890, see
\cite{KleinGesammelteIII}, p.~577.} be called {\em automorphic
functions}. Because of excessive overwork and asthma attacks
Klein suffered a breakdown of health in 1882, from which he only
slowly recovered. After that, he felt unable to continue his
research as before,\footnote{See Klein's own recollections in
\cite{KleinGesammelteII}, p.~258 and \cite{KleinGesammelteIII},
p.~585.} and decided to concentrate on promoting the subject in
lectures and textbooks, thus trying to attract young
mathematicians into this field.

During that time Klein wrote his {\em Vorlesungen \"uber das
Ikosaeder} \cite{Ikosaeder}, which was published in 1884. This
book was to be the first of a multivolume project on ``regular
solids, modular functions, and automorphic
functions''\footnote{See \cite{FrickeKleinI}, p.~V of the {\em
Vorrede}, and \cite{KleinGesammelteIII}, p.~742.}. In the autumn
of 1887 Klein invited Fricke to join these efforts. While
Fricke's prowess as a mathematician was already apparent in his
thesis of 1885, this project would be the stepping stone for his
academic career. Klein and Fricke went on to write a two-volume
book on {\em Elliptic Modular Functions}
\cite{KleinFrickeI,KleinFrickeII} published in 1890 and 1892, the
first of which Fricke used to fulfill the requirements for his
{\em Habilitation} in Kiel. These were followed in 1897 by {\em
Automorphic Functions~I} \cite{FrickeKleinI}.

At the start of their cooperation, Klein, as a well established
expert, was the main author, while Fricke as a young Ph.D.\ and
aspiring mathematician meticulously worked out the details. Over
the years, as Fricke gained experience and his scientific
reputation grew, the character of their collaboration changed.
While Fricke more and more exerted influence on contents and
style, thus in fact becoming the main, and in large part the only
author and the driving force, Klein remained as the initiator and
public figurehead of the project and acted as some kind of
spiritual guide in the background, every once in a while slightly
changing the direction or contents.\footnote{See
\cite{FrickeKleinII}, p.~VI et seq. of the {\em Vorrede}.} This
change of role was reflected by the change of order in which the
names of these two authors appeared on the title page of {\em
Automorphic Functions~I}.

\section{William Burnside, and automorphic functions}

The life of William Burnside is fairly well documented in \cite{BurnsideOb,BurnsideCollected1},
so we may limit ourselves to a short survey.

Burnside was born on 2 July 1852 in London where he also grew up
and went to school. In October 1871 he started studying
mathematics in Cambridge. He graduated from there in 1875, and
became a fellow of Pembroke College. In 1885 Burnside was
appointed professor of mathematics at the Royal Naval College at
Greenwich, a position which he held until his retirement in 1919.
He continued his research until his death in 1927.

Burnside's preoccupation with automorphic functions began about
1891. In that and the following year he authored a sequence of
three papers
\cite{BurnsideDisco,BurnsideAutomorpheI,BurnsideAutomorpheII} on
this subject. These marked a turning point in Burnside's
development as a mathematician. Starting out in the
Cambridge tradition of applied mathematics, with a personal
emphasis laid on kinematics, kinetics, and especially
hydrodynamics, he had now devoted large parts of these articles,
which were going to exceed all of his other research publications in length,
to a more formal approach. Furthermore the papers stood at
the beginning of his
growing involvement in the theory of groups, mainly in discrete
and finite groups.\footnote{This estimation was already given in
Burnside's obituary \cite{BurnsideOb}, p.~76ff.}

Both to elaborate the connection between automorphic functions and group theory,
and to assess Burnside's achievements, at this point it is advisable to give
a rough outline of the mathematics involved:

Let $\Gamma$ be a discrete group of linear fractional transformations
$z\mapsto \frac{\alpha z + \beta}{\gamma z +\delta}$ with
$\alpha,\beta,\gamma,\delta\in\mathbf{C}$ and $\alpha \delta - \beta\gamma = 1$,
the topology being induced by the Euclidean metric.

If $d$ is an arbitrarily chosen integer, an {\em automorphic form} of dimension $d$
with respect to $\Gamma$ is a meromorphic function $\theta$ defined on a domain
in $\mathbf{C}$ satisfying the functional equation
$$
\theta\left(\frac{\alpha z +\beta}{\gamma z + \delta}\right) = (\gamma z + \delta)^d\cdot \theta(z)
$$
for all elements $z\mapsto\frac{\alpha z+\beta}{\gamma z+\delta}$ of $\Gamma$,
whenever $\theta$ is defined.

An automorphic form of dimension 0 is called an {\em automorphic function}.
Thus, automorphic functions with respect to $\Gamma$ are invariant under the group
action of $\Gamma$.

Given a rational function $z\mapsto H(z)$ and a negative integer $d$,
the {\em Poincar\'e series} of dimension $d$ (with respect to $\Gamma$ and $H$) is defined as
$$
\sum_k H\left(\frac{\alpha_k z+\beta_k}{\gamma_k z+\delta_k}\right)\cdot(\gamma_k z+\delta_k)^d,
$$
where the summation is taken over all of the  elements
$z\mapsto\frac{\alpha_k z+\beta_k}{\gamma_k z+\delta_k}$ in the
group $\Gamma$ (because $\Gamma$ is discrete the summation is
taken over at most countably many elements).

Whenever such a series is uniformly convergent in an open domain,
it represents an auto\-morphic form of dimension $d$ in that domain.
Therefore the question arises quite naturally of how the
convergence of Poincar\'e series depends on the group $\Gamma$
and the dimension $d$. Poincar\'e himself
\cite{PoincareActa1,PoincareActa3} had already shown that for
dimensions $d\le -4$ these series are uniformly convergent where
defined.

In \cite{BurnsideDisco,BurnsideAutomorpheI,BurnsideAutomorpheII}
Burnside significantly extended Poincar\'e's results. Burnside
started from a concrete physical problem in hydrodynamics: the
uniform streaming motion in two dimensions around a set of
circles. The physical description of uniform flow together with
the boundary conditions given by the circles leads to a
differential equation the solution of which might be represented
by a Poincar\'e series. Burnside subsequently obtained results on
the convergence of Poincar\'e series of dimension $d=-2$ which
also extend to $d=-1$.

Although with these papers he entered a field of research which
was new to him, Burnside took notice of other important
publications in this direction besides Poincar\'e's. During the
course of his articles he incorporated Klein's freshly created
terminology of {\em automorphic function} and {\em Primform}.
Furthermore he mentioned a paper \cite{SchottkyPoincare} by
Friedrich \mbox{Schottky} (1851--1935), which had advanced the
theory of Poincar\'e series considerably. Nevertheless
Burnside's methods differ substantially from Schottky's, and thus
lead to parallel and in some instances improved results.

\section{Fricke composing {\em Automorphic Functions II}}

In the autumn of 1899, the year after his and Burnside's article on the simple group of order 504,
Fricke focussed on working out the structure and the details of {\em Automorphic Functions~II}.
This volume was planned to be the completion of his and Klein's project.
At the end of that year, Fricke had written the first one and a half chapters,
approximately 100 pages, and had sketched an outline of the third,
which was intended to cover the theory of Poincar\'e
series.\footnote{Fricke to Klein, 9 Sept, 8 Nov, 27 Dec 1899, Klein Nachlass 9, SUB G\"ottingen.}
For this chapter, Fricke collected the published convergence results known to him,
starting from Poincar\'e's fundamental papers \cite{PoincareActa1,PoincareActa3}.
For dimension $d=-2,$ he took into consideration Schottky's already mentioned article
\cite{SchottkyPoincare} from 1887, and the Ph.D.~thesis \cite{RitterDiss} of Ernst Ritter
(1867--1895) from 1892, which had been written in G\"ottingen under the auspices of Klein.
By extending the ideas which he found in these sources, mainly using Schottky's method of proof,
Fricke managed to fill several gaps by the summer of 1900.
His progress and the development of the book can be traced almost minutely through a series
of letters to Klein.\footnote{Fricke to Klein, 25 Apr, 12, 19, 25 May, 1, 8, 17 July, 11 Oct 1900,
Klein Nachlass 9, SUB G\"ottingen.}

Although Burnside's results from 1891 and 1892 were highly relevant to the subject, they are not
mentioned in Fricke's correspondence, and seem to not have been part of the initial manuscript.
Fricke did report to Klein on an article of Whittaker on automorphic
functions,\footnote{Fricke to Klein, 23 Dec 1898, Klein Nachlass 9, SUB G\"ottingen.}
in which Burnside's achievements were explicitly mentioned,
but he only referred to a preliminary abstract in~\cite{WhittakerAutomorpheAbstract}.
If he had considered the later article \cite{WhittakerAutomorphe} from 1899,
he certainly would have come across Burnside's papers.

\medskip
While streamlining the exposition of his manuscript, Fricke had
also excerpted two short articles
\cite{FrickeRitter1,FrickeRitter2}, which he published in {\em
G\"ottinger Nachrichten}. So exultant was he about his
accomplishments that on 17 July 1900, after months of intensive
work, he wrote to Klein:\footnote{Fricke to Klein, 17 July 1900,
Klein Nachlass 9, SUB G\"ottingen (our translation).} ``Only now
have I gained confidence that `Automorphic Functions II, Part~I'
will actually succeed.''

Despite this show of optimism, Fricke and Klein both lived in
constant trepidation that the general interests of the
mathematical community had gradually but surely shifted away from
automorphic functions towards other, more foundational
issues.\footnote{Only a month after the enthusiastic letter from
Fricke, David Hilbert (1862--1943) delivered his famous talk {\em
Mathematische Probleme} before the Second International Congress
of Mathematicians in Paris, which included a list of open problems
that would subsequently be considered as the most important to
``Modern Mathematics''.} Therefore Fricke planned to present his
new results on Poincar\'e series to the general mathematical
public at the meeting of the {\em Deutsche
Mathematiker-Vereinigung}, which was to take place in Aachen
from 16 to 23 September 1900. He regarded this talk as an
opportunity both for promoting his work and for probing the
possible reaction towards the forthcoming book on automorphic
functions.

However, when he returned from Aachen, Fricke was in a state of turmoil.
He wrote to Klein:\footnote{Fricke to Klein, 11 Oct 1900, Klein Nachlass 9,
SUB G\"ottingen (our translation).}
``This year's meeting did not turn out very propitious for me, as I
did not find the right kind of support for the areas of interest
with which I came. I came there `automorphically charged', so to
speak, but in Aachen, as well as now afterwards, I was left with
the impression that I will commit an anachronism with the second
volume on automorphic functions.''

Fricke further described how he reacted to this blow to his confidence and which actions he took:
``It is nevertheless a case of `the sooner the better' and I would have to wage a pretty fierce battle
if I intended to definitively abandon this project.
My policy therefore was to make myself as free as possible from other hindrances
so as to develop auto\-morphic functions more quickly.
I have therefore asked Lampe\footnote{Emil Lampe (1840--1918) was the editor of
{\em Jahrbuch \"uber die Fortschritte der Mathematik}, the main reviewing journal in mathematics
at that time, which Fricke voluntarily worked for as a reviewer of articles in number theory.}
to nominate another colleague as reviewer for number theory starting next year [$\ldots$]
On the other hand, following my mood I~have begun every possible
correspondence concerning automorphic functions.''

Thus in defiance of both his own disenchantment and the seeming
lack of interest of the mathematical community, Fricke
concentrated more than ever  on concluding the first part of {\em
Automorphic Functions II}, which at that time he only had to put
the finishing touches to. \mbox{After} returning from a short
journey to Paris, he sent letters to several experts on
auto\-morphic functions, among them David Hilbert\footnote{Fricke
to Hilbert, 30 Sept 1900, Hilbert Nachlass 107, SUB
G\"ottingen.}, Ferdinand Lindemann\footnote{Lindemann to Fricke,
22 Nov 1900, Fricke Nachlass, UA Braunschweig.} (1852--1939),
Joseph Wellstein\footnote{Wellstein to Fricke, 29 Oct 1900, Fricke
Nachlass, UA Braunschweig.} (1869--1919), and Wilhelm
Wirtinger\footnote{Wirtinger to Fricke, 8, 15 Oct 1900, Fricke
Nachlass, UA Braunschweig.} (1865--1945). He had met all but
Lindemann in Aachen, and in each case the correspondence was used
to continue discussions from the meeting.

\section{A letter of response and enquiry}

It seems that until the Aachen meeting Fricke was completely
unaware of Burnside's work on automorphic functions. As a
mathematician, Burnside was known to him at least since the
episode of their nearly simultaneous working on the simple group
of order 504 of two years before. In fact, Fricke had been chosen
by Klein during the process of editing Burnside's paper to
evaluate its correctness, and thus had been spurred into
activity, which had resulted in his own article in {\em
Mathematische Annalen} shortly afterwards.\footnote{Fricke to
Klein, 1, 12, 13 Sept, 13 Oct 1898, Klein Nachlass 9, SUB
G\"ottingen.}

After Aachen, however, Fricke at last was aware of Burnside's
results, as can be seen by the summary of his Aachen talk in
\cite{FrickePoincareDMV}, where he clearly states them. We may
assume that he was informed about Burnside's articles during the
meeting by Heinrich Burkhardt (1861--1914), a close associate of
Fricke and Klein. Burkhardt had reviewed Burnside's papers for
{\em Jahrbuch \"uber die Fortschritte der Mathematik} in 1895,
which was immediately followed by another review from him
discussing Ritter's thesis.\footnote{See {\em Jahrbuch \"uber die
Fortschritte der Mathematik}, Volume 24, p.~391ff.} Presumably
this incident in Aachen might have caused Fricke to start his
correspondence in order to avoid other such embarrassing
situations due to ignorance of relevant literature.

Thus, when in short succession he sent off the series of letters
mentioned above, he also started an exchange of ideas with
Burnside. Unfortunately, of this particular correspondence, only
two of Burnside's letters seem to have survived. The first one of
these\footnote{Burnside to Fricke, 28 Oct 1900, Fricke Nachlass,
UA Braunschweig.} starts as follows:

\medskip
\begin{quote}
\begin{em}
  \hfill                        Bromley Road\\
\phantom{m}  \hfill                             Catford\\
\phantom{m}  \hfill                             Oct 28, 1900.\\
\noindent
Dear Professor Fricke

 Many thanks for your last
very kind letter. I hope very much you will
carry out your intention of making a trip to
England; and if you do, I trust you will let
me shew you some hospitality. It would give me
great pleasure if you would come and stay here for
a few days so that we might become personally
acquainted. I can read German with comparative
ease, but I cannot speak it at all {\sl\&} I fear now
it is rather late to learn.

 At the risk of wasting your time I must return
once more to the subject of my last letter, to free
myself from the reproach of having forgotten the
main point of my argument. Since I wrote you
I have looked up Poincar\'e's article in Acta Math.
Vol~I, with which I was more familiar in 1892 than
I am now; and I have reconstructed my argument.
I put it forward without much confidence, but at
all events it is what I~wanted to send you in
my last letter.

 Very many thanks for the papers from the
G\"ottingen Nachrichten you have sent me. I am
posting with this a short paper on groups of
finite order which I contributed to the Stokes
memorial volume.
\phantom{m}\hfill       Yours very sincerely\\
\phantom{m}\hfill               W. Burnside.

\end{em}
\end{quote}

\medskip
It is safe to assume that the papers from {\em G\"ottinger
Nachrichten} which Fricke sent to Burnside are those articles
\cite{FrickeRitter1,FrickeRitter2} mentioned above which Fricke
published while working out the details of {\em Automorphic
Functions II}.

Returning Fricke's kindness, Burnside sent to him his contribution \cite{BurnsideStokes}
to volume 18 of {\em Transactions of the Cambridge Philosophical Society}
which had been reserved for articles dedicated to George Gabriel Stokes (1819--1903)
on the occasion of the fiftieth anniversary of his being Lucasian Professor of Mathematics
at the University of Cambridge in October 1899.

The reconstructed argument that Burnside mentioned
constitutes the contents of the second page of Burnside's letter:

\begin{quote}
\begin{em}
Putting
$$f(z)=\frac{\sum N(z_i)(\gamma_i z+\delta_i)^{-2m-2}}{\sum D(z_i)(\gamma_i z+\delta_i)^{-2m}\enspace\enspace}$$
take $\sum D(z_i)(\gamma_i z+\delta_i)^{-2m}$ to be free from poles
then p.~227, it has $2m(n-1)$ zeros.

 Again taking $\sum N(z_i)(\gamma_i z+\delta_i)^{-2m-2}$ free of poles
it can (p.~282) be expressed linearly in
terms of $(2m+1)(n-1)$ linearly independent
functions: i.e.\ the numerator is a linear
homogeneous function of $(2m+1)(n-1)$ arbitrary
constants. If each zero of the den[ominator] is
a zero of the numerator there are $2m(n-1)$
linear relations among the constants,
leaving over $n-1$ arbitraries. But from
the physical considerations $f(z)$ ought to
contain $n$ arbitraries, viz.\ the $n$
circulation constants. Hence I inferred
that the function postulated physically
cannot be represented in the form

$$\frac{\sum N(z_i)(\gamma_i z+\delta_i)^{-2m-2}}{\sum D(z_i)(\gamma_i z+\delta_i)^{-2m}\enspace\enspace}.$$

\smallskip
{\rm [On the margin of this page we find the remark:]}

references to Poincar\'e's paper in
 Vol I of Acta Mathematica.
\end{em}
\end{quote}

\medskip
Let us try to put this proof in the context of Burnside's original
paper. Obviously, Fricke's question pertained to certain
discontinuous groups, the fundamental domains of which are
bounded by $n>1$ pairs of circles. In the last parts of his paper
\cite{BurnsideAutomorpheI}, p.~287 et seq., Burnside had asserted
by counting dimensions that there exist everywhere finite
automorphic forms of dimension $-2$ which cannot be represented
by Poincar\'e series. Fricke seems to have asked if such an
automorphic form could be represented by a quotient of Poincar\'e
series of appropriate dimensions, thus avoiding any difficult
convergence considerations. Burnside replied that such a
representation is not possible in general, giving as a reason an
extended version of his dimension argument. Since Fricke at that
time still collected any information concerning automorphic
functions in connection with Poincar\'e series, he was clearly
interested in Burnside's argument.\footnote{Fricke himself had
explicitly expressed an analogous assertion in one of his earlier
major papers in 1892 (see \cite{Fricke237_247}, p.~453), where he
referred to \cite{RitterDiss}. Remarkably, Fricke would later use
the findings of that paper in his research on the simple group of
order 504.}

These pages of Burnside's letter enable us to restore the basics of the correspondence so far:
Fricke started by asking a mathematical question on the representability
of automorphic functions by Poincar\'e series. Burnside replied by stating the general
impossibility, admitting that he could not remember the details of his argument.
Fricke answered with a ``very kind'' letter in which he included personal details,
such as his plan to travel to England.
Maybe he expressed his wish to get personally acquainted with John Perry (1850--1920),
whose {\em Calculus for Engineers} he was concurrently translating into German.
Perhaps he mentioned as well that he had been invited to Cambridge in June 1899
to join the celebration in honour of Stokes, but had not attended.
He also included his latest scientific publications.
Burnside felt himself obliged to look up the missing details of his proof,
and supplemented his reply with one of his latest publications.

\medskip
But Burnside's letter does not end with his more comprehensive
answer to Fricke's question, as would usually be the case in a
correspondence between two mathematicians who barely know each
other. Instead, Burnside added a third page which reads as
follows:

\medskip
\begin{quote}
\begin{em}
 I take the opportunity of asking you, whether the
following question has ever present\-ed itself to
you; and if it has, whether you have come to
any conclusion about it.
\smallskip

 Can a group, generated by a \underline{finite} number
of operations, and such that the order of
every one of its operations is finite and
less than an assigned integer, consist of
an infinite number of operations.

 E.g.\ as a very particular case:--

If $S_1$ and $S_2$ are operations, and if
$\Sigma$, representing in turn any and every combination or repetition
of $S_1$ and $S_2$, such as $S_1^a S_2^b S_1^c\ldots S_2^e$, is
such that $\Sigma^m=1$, where $m$ is a given
integer, is the group generated by $S_1$ and
$S_2$ a group of finite order or not.
Of course if $m$ is $2$, the group is of order $4$
and if $m$ is $3$ the group is of order $27$;
but for values of $m$ greater than $3$, the
question seems to me to present
serious difficulties however one looks at it.
\end{em}
\end{quote}
\medskip

In a more specified form, the question raised above would later
be called the {\em Burnside Problem}. Although it is presented
here in statu nascendi and for the special case of two generators,
the fundamental problem is already clearly stated.

\section{A letter of thanks and report on progress}

After Fricke and Klein had finished their manuscript in December 1900,
the first part of {\em Automorphic Functions II} was published in May 1901.
As was customary at that time, they sent off a number of copies of their book
to colleagues who worked in that particular area of mathematics
or to whom they wanted to express some measure of gratitude or friendship.
Due to the high costs, the list of people thus recognised had to be
restricted in number to a select few.
On that specific occasion, Burnside belonged to the recipients of such a copy.
Usually, this honour was answered with short polite letters of thanks,
the tone of these missives determined by the degree of personal acquaintance to the author.
Again, Burnside's reply\footnote{Burnside to Fricke, 9 June 1901, Fricke Nachlass, UA Braunschweig.}
goes beyond this:

\medskip
\begin{quote}
\begin{em}
\hfill    The Croft,\\
\phantom{m}\hfill       Bromley Road,\\
\phantom{m}\hfill               Catford.\\
\phantom{m}\hfill            June 9, 1901.\\[0.5ex]
\noindent
Dear Prof. Fricke

 Very many thanks
to you and to Prof.\
Klein for your kindness
in sending me a
copy of the new
instalment of your
book on automorphic
functions. I got the
first volume when it
appeared and have
read the greater part
of it; and I look
forward with great
pleasure to studying
the new part.

 I have recently returned
to the question I wrote
you about in the
winter, viz.\ that of the
discontinuous group
defined by
     $$S^m=1$$
when $S$ represents any and
every combination of $n$
independent generating operations
 $A_1,A_2,\ldots, A_n$;
and $m$, $n$ are given integers.

I find that when
    $m=3$
and $n$ is given, the order
can be determined by a
kind of recurring formula.
In particular if $n=3$
the order is $3^{17}$.

 For $m=4$, $n=2\;$ I find
$2^{12}$ for the order. These
results, if correct as I
believe them to be, would
seem to shew that a
graphical method, i.e. the
consideration of the
network of polygons, is
practically out of the
question from the extreme
complexity of the figure with
so great a number of
polygons.

 So far I am quite
baffled by the case
$n=2$, $m=p\,$, a prime greater
than $3$; but it is easy
to shew that the order
cannot be less than $p^{p+2}$
and I think it is
probably greater, if finite
at all.\\
\phantom{m}\hfill          Believe me\\[0.5ex]
\phantom{m}\hfill        Yours very sincerely\\
\phantom{m}\hfill          W. Burnside.
\end{em}
\end{quote}

\medskip
With the letters from Burnside at hand, we can try to reconstruct
Burnside's progress and the genesis of his famous paper
\cite{BurnsideProblem} of 1902 where Burnside's results were
published in a coherent and detailed form. Below we let $B(m,n)$
denote the universal group with $m$ generators, all elements of
which have an order dividing $n$, which today is appreciatively
called the {\em Burnside group} of exponent~$n$ (with
$m$~generators).

It seems that in October 1900 Burnside had not progressed very far
in determining answers to his own question. In the first letter
Burnside presented the two examples $B(2,2)$ and $B(2,3),$ the
first of which is quite obvious while the other one was known to
him for several years: $B(2,3)$ had already appeared as an
illustrative example of the graphical method of analysing groups
in his article \cite{BurnsideGraphic} in 1893, and was used again
in the same context in his monograph on finite groups
\cite{BurnsideGroups} in 1897. This method, which Burnside
described quite thoroughly in his book, starts from an abstract
finitely presented group and leads to a representation of it as a
group of linear fractional transformations acting on the
plane.\footnote{Burnside attributed the origin of this method to
the fundamental paper {\em Gruppentheoretische Studien}
\cite{DyckI} of Walther (von) Dyck (1856--1934). Dyck had been a
student and an assistant of Klein in Munich and in Leipzig, where
he wrote \cite{DyckI} as his {\em Habilitation} thesis. The title
of nobility ``von'' was bestowed on Dyck in 1901.} Burnside's
1893 paper~\cite{BurnsideGraphic}, which was also his first
dedicated solely to group theory, may be interpreted as
containing the seeds from which the Burnside Problem
originated.\footnote{In \cite{BurnsideGraphic} Burnside
considered the infinite group given by generators $P$ and $Q$
which satisfy
$$P^3=Q^3=(PQ)^3=E$$
and asked which further relations enforce the resulting group to be finite.}

When Burnside wrote to Fricke again in June 1901, he could report
on some substantial advance, and he was even quite close to the
final article. Having recapitulated the main problem in the
special form that can also be found in the article (though by
1902 the role of the integer parameters $n$ and $m$ would be
interchanged), Burnside continued to give an account of his
results up to that time, leaving out any proofs: First of all
Burnside told of a ``recurring formula'' for the order of
$B(m,3)$. If we presume that the erroneous value $3^{17}$ was
inadvertently written in the letter (and not simply due to an ink
spot), we contend that he intended to give as a particular
example $|B(3,3)|= 3^7$. While the recurrence formula which
Burnside presented in his later article leads to the correct
order of $B(3,3)$, it still only gives an upper bound for greater
numbers of generators, as was shown by later authors in
\cite{LeviWaerden}. The next result Burnside stated in his
letter, $|B(2,4)|=2^{12}$, was right on the spot (although we now
know that the proof in \cite{BurnsideProblem} only suffices to
show that the order of this group divides $2^{12}$). Finally
Burnside wrote of a lower bound for the order of $B(2,p)$, $p$
being a prime greater than 3, always keeping in mind that the
order might even be infinite. Using facts in elementary group
theory, he was able to improve on this result by the time of the
final version.

It seems that, while Burnside was on the right track in all cases,
sometimes he was a little too confident and stepped into hidden
pitfalls. While the methods and proofs presented in the final
paper were correct, they did not suffice to show all parts of the
stated claims. With hindsight, this does not diminish Burnside's
achievements. As Burnside told Fricke in his second letter, the
orders of the groups appearing force the researchers to go beyond
the scope of any visually oriented method. Thus here  a purely
abstract approach was truly necessary, which Burnside provided
splendidly.\footnote{Even the proof readers of the final article
seemed to be slightly overwhelmed by its content, since they
missed Burnside's accidental switch in the name of generators
from $A,B$ to $P,Q$ in the paragraph before last.}

\medskip
Both in the letters and in the article \cite{BurnsideProblem} Burnside treated
special instances of the question which today is justifiably called the {\em Burnside Problem}:
Is a group with a finite number $m$ of generators necessarily finite
if the orders of its elements all divide a given integer $n$?

However, the introductory phrase of Burnside's first preserved letter addresses
a related, but slightly different question: Is a finitely generated group finite
if the orders of its elements are less than an integer $n$?
This particular variant did not see print before the publication of the textbook
of Harold Hilton (1876--1974) on finite groups \cite{HiltonGroups} in 1908.
There it could be found in the appendix, where it was part of a list of twelve research problems,
most of which Hilton expressly attributed to Burnside.\footnote{The first to
address this problem was B.~H.~Neumann in 1937 \cite{NeumannBounded},
who termed it the {\em bounded order problem}.}

Moreover when in the second letter Burnside repeated the problem as a small reminder,
he explicitly put it within the context of discontinuous groups, just as in the final article.
For somebody belonging to a later generation of mathematicians this seems to be quite unnecessary,
because the Burnside Problem can be seen as a problem in abstract group theory,
totally divorced from any geometric interpretation.

Strangely enough, in the initial sentence of his 1902 article,
Burnside kept his reference to discontinuous groups and instead
dropped both the assumption of existence of a universal upper
bound for the orders of each group element, and, more grievously,
the prerequisite that the group in question be finitely
generated. The last omission is substantial, because obviously
any infinitely generated abstract group is infinite, even if all
element orders are finite. Abandoning only the first condition
leads to a generalisation of the original problem, which nowadays
is called the {\em General Burnside Problem}: Is a finitely
generated group necessarily finite if the orders of its elements
are finite?

The vagueness of the first sentence of \cite{BurnsideProblem} has given rise to a number of
interpretations including a discussion if Burnside's reference to discontinuous groups
already implied those essential conditions.\footnote{See e.g.\ \cite{NewmanBurnside}.
In the case of complex matrix groups, Burnside himself has shown in 1905 \cite{BurnsideLinearFinite}
that the bounded order problem has a positive solution without resorting
to the further assumption of the group being finitely generated.}
To shed some more light on this particular reference, we have to take a closer look
on Burnside's understanding of groups.

\section{Burnside's attitude towards groups -- and towards Fricke}
Nearly all of Burnside's publications suggest that he thought of
groups as being composed of ``operations''. This observation has
already been pointed out in case of his book by Neumann in
\cite{NeumannBurnside}. There it was contrasted with the
definition of a group in the second volume of {\em Lehrbuch der
Algebra} \cite{WeberAlgebraEd1} of Heinrich Weber (1842--1913),
where the author more abstractly only speaks of ``elements'' as
constituents of a group.

When at the turn of the century Burnside's main interests shifted
to abstract group theory, his approach to (countably) infinite
groups always tended to be in connection with what at that time
were called ``discontinuous groups''. The contemporary
understanding of this term differed from today's: it was mainly
used to discriminate this class of groups from ``continuous
groups'', which at that time more or less meant Lie groups.
Burnside tried to give a precise definition of both concepts when
he was given the opportunity of writing an expository article on
the theory of groups in the {\em Encyclopaedia Britannica}
\cite{BurnsideEncBrit10,BurnsideEncBrit11}.

According to our analysis his 1893 paper \cite{BurnsideGraphic}
seems to be his only research paper where Burnside speaks of the
elements of a group as ``symbols'' and not as ``operations''. Its
central problem is indeed formulated as an abstract word problem.
However, the solution is deduced in large part via the graphical
method, which is directly connected to discontinuous groups. In
his 1901 letter to Fricke we find the explicit statement about
the unsuitability of this approach, which again confirms our
conclusion that \cite{BurnsideGraphic} can be seen as a direct
precursor of \cite{BurnsideProblem}.

Although Burnside was well aware of the difference between
abstract finitely generated groups and concrete discontinuous
groups, most of the time he did not seem to feel the necessity to
make a clear-cut distinction between these two classes. Both the
evidence of the letters to Fricke and his 1893 paper show that
Burnside's interest in the Burnside Problem originally stemmed
from concretely given discontinuous groups. Thus, while today his
1902 paper \cite{BurnsideProblem}, where he first presented his
set of problems to a general public, for the most part is read as
a piece of pure abstract group theory, he himself explicitly put
his ``unsettled question'' into the context of discontinuous
groups. We assume that he did so, as he had done several months
before in his letter, without giving a second thought to this
wording.

\medskip
These insights into Burnside's view about groups enable us to
explain why he had used the occasion of answering a mathematical
question on automorphic functions to present his current research
problems in the field of group theory, together with the results
he had achieved, to a virtual stranger such as Robert Fricke. At
first sight, his behaviour seems surprising, because, as we noted
above, he was known to have worked more or less in seclusion,
without many personal contacts with other non-British
mathematicians, and his contact with Fricke happened rather by
chance than by his own choice.

Although he did not know Fricke personally, Burnside was well
acquainted with the subjects and results of Fricke's work.
Burnside was said to be ``well informed about the published
literature''\footnote{See \cite{NeumannBurnside}, p.~32.}, and
he himself wrote in his letter that he could ``read German with
comparative ease''. We can thus safely assume that the remark in
his letter that he had read the greater part of {\em Automorphic
Functions I} was not just a set phrase, but a fact. Indeed, the
circle of ideas encompassing automorphic functions, discrete
groups and physical applications thereof was not only well known
to him, but for some time in his life formed one of his central
research interests. We agree with Forsyth's
assessment\footnote{See \cite{BurnsideOb}, p.~77.} that Burnside
initially encountered infinite discrete groups in this context.

Moreover, Burnside's opinion of how expository mathematical writing should be done
was quite close to Klein's and Fricke's illustrative approach.
On the occasion of a presidential address \cite{BurnsideAddress}
to the London Mathematical Society in 1908, Burnside characterized
the succinct way of exposition in research papers on abstract group theory
as ``driest formalism'', made up of ``a series of conundrums''.
If his own article \cite{Burnside504} on the simple group of order~504
were assessed according to the criteria voiced in the address,
only the last paragraph, where Burnside made the connection with a concrete example,
would save this paper from being
``a merely curious illustration of non-commutative multiplication''.

On the other hand, Burnside's presidential address also describes the way in which
such an article should be written: ``The reader is led naturally from the concrete
to the abstract, and acquires by actual instances the new ideas necessary to the theory.''
Fricke's article \cite{Fricke504} on the simple group of order~504 exemplifies
to the utmost Burnside's demand for mooring abstract results in concrete examples,
and appears as if written in full accordance with Burnside's guidelines.

When Burnside himself wrote on the theory of groups for the 10th edition
of {\em Encyclopaedia Britannica} in 1902, he quoted as authorities
on both discontinuous and finite groups Felix Klein with his book
on the icosahedron as well as several others of Klein's school:
Dyck with his {\em Gruppentheoretische Studien}, Josef Gierster (1854--1893),
and, final\-ly, Fricke with {\em Elliptic Modular Functions I}
and {\em Automorphic Functions I}.
Later on, in the next edition of {\em Encyclopaedia Britannica} in 1911,
while severely cutting down the number of references,
Burnside left the books of Klein and Fricke on his list,
and added {\em Automorphic Functions~II}.

In summary, Burnside's personal development towards group theory´
paralleled that of a number of Klein's former students.
The mathematicians affiliated to this group were the natural choice
as addressees for Burnside's questions.
When the occasion arose and Burnside came into direct contact with
one of its more exposed members, namely Robert Fricke,
who had shown a profound versatility in those areas close to his heart,
he seized the opportunity and allowed a deeper insight into his own current research
with high hopes of having found a like-minded spirit and highly competent partner.

\medskip
As a matter of fact, Burnside's inquiry was  to produce a prompt
reaction from Fricke, as is shown by the reverse sides of the
pages of the letter from 1900. While they initially had been left
blank by Burnside, Fricke had used them as some kind of scratch
paper to sketch his immediate ideas on Burnside's question. Thus
on the reverse of the third page we find a pencil sketch of a
regular hexagon consisting of $54$ alternating black and white
equilateral triangles, which forms the graphical representation
of the group $B(2,3)$ and which can also be found as part of the
only accompanying figure of Burnside's 1893 paper
\cite{BurnsideGraphic}. We are safe to assume that Fricke did not
need to look it up but was able to produce it by himself.
Furthermore on the reverse of the second page there are several
computations in Fricke's handwriting. Obviously they were done to
work out a hint on the solution of the first open case, the group
of exponent~$4$ generated by two elements, by calculating all
$2\times 2$ matrices with entries in the integers modulo~8
generated by $1\;2\choose0\;1$ and $1\;0\choose2\;1$. Since this
did not lead to any useful insights, Fricke eventually abandoned
this line of research.

\medskip
The remarkable similarity between Burnside's and Fricke's grasp of
the concept of a group became visible again several years later.
When after Weber's death his {\em Lehrbuch der Algebra} went out
of print, the publisher offered Fricke the opportunity to
elaborate a new revision with the same title, based on the
original. At the point where Fricke gave the definition of a
group, he closely followed Weber, but added the following -- from
a mathematical point of view unnecessary --
introduction\footnote{See \cite{FrickeAlgebraI}, p.~267 (our
translation and emphasis).}, which sounds only too reminiscent of
Burnside:
``Let there be given a system of like {\em operations}, or like analytic expressions,
or mathematical entities described in any other way,
which to refrain from any interpretation we call `elements' [$\ldots$]''

\section{Aftermath}

In {\em Automorphic Functions II}, Fricke explicitly acknowledges
Burnside's contribution to the convergence of Poincar\'e series
of dimension $-2,$ but the results are only mentioned in a
reporting tone and seem to have been worked into the manuscript
of the book in a very late phase.\footnote{See
\cite{FrickeKleinII}, pp.~157, 160, 166.} The proof presented
there follows the lines of Schottky. However, after finishing his
manuscript, Fricke seems to have worked out, simplified, and
generalized Burnside's results to prove the convergence of
certain Poincar\'e series of dimension $-1$. How important Fricke
deemed these results to be can be seen by the fact that he took
the chance to include these achievements in the {\em Festschrift}
\cite{FrickePoincare} which was published on the occasion of
Richard Dedekind's 70th birthday in October 1901.

The publication of the further parts of {\em Automorphic Functions II} was considerably delayed:
Fricke and Klein wanted to include some key theorems, proofs of which were
independently given by Poincar\'e and Klein in the early 1880s.
Unfortunately, by the turn of the century, these proofs were not considered
sufficiently rigorous anymore.\footnote{It seems that Ernst Zermelo (1871--1953)
had pointed out that there were major gaps in the initial proofs.
See Fricke to Klein, 9 Sept 1899, 1 Sept 1901, Klein Nachlass 9, SUB G\"ottingen.}
Only after 1907, when L.~E.~J.~Brouwer (1881--1966) and Paul Koebe (1882--1945)
had achieved new results on the uniformization of algebraic curves,
based on the advances of modern point-set topology,
could Fricke continue completing {\em Automorphic Functions II}.
The final two parts were published in 1911 and 1912, respectively.

After their exchange of letters, despite their similar attitudes,
Burnside and Fricke seem to have lost contact with each other.
Although Fricke travelled to England in spring 1903 to meet John
Perry, we have found no evidence that he took up Burnside's
invitation. On the other hand, Burnside became a member of the
{\em Deutsche Mathematiker-Vereinigung} in 1904, and remained so
for the rest of his life, even through the First World War, thus
always keeping in touch with the mathematical community in
Germany.

Both the Burnside Problem and the problem of convergence of Poincar\'e series
were the source of inspiration for later mathematicians.
In both cases there are still open questions awaiting their definitive solutions.

\section{Acknowledgements}
The letters referred to in the text may be found either in the university archives (UA)
in Braunschweig, or in the university library (SUB) in G\"ottingen, respectively.
The authors are indebted to several staff members of these institutions
for their advice, support, and patience.

Furthermore they would like to thank Jeremy Gray, Tony Mann, and Peter M.~Neumann
who read various preliminary versions of this article and were so kind
to help weeding out some mistakes.
Clearly the authors bear the responsibility for any remaining errors.

Finally the second author wants to acknowledge the financial support
by Johannes Fuhs without which this research would not have been possible.

\nocite{BurnsideCollected1}
\nocite{BurnsideCollected2}
\newcommand{\noopsort}[1]{}
\bibliographystyle{amsplain}

\end{document}